\documentclass[final]{dmtcs-episciences}

\usepackage{amsmath}
\usepackage{amssymb,amsfonts,amsthm}
\usepackage{tikz,graphicx}

\usepackage[utf8]{inputenc}
\usepackage{subfigure}
\newtheorem{theorem}{Theorem}[section]
\newtheorem{corollary}[theorem]{Corollary}

\newtheorem{proposition}[theorem]{Proposition}
\newtheorem{conjecture}[theorem]{Conjecture}
\newtheorem{remark}[theorem]{Remark}
\newtheorem{observation}[theorem]{Observation}
\newtheorem{claim}{Claim}[theorem]

\newcommand{\w}{{\rm w}}
\newcommand{\cF}{{\cal F}}

\newcommand{\residual}{G_S^\iota}

\newcommand{\resT}{T_S^\iota}

\newcommand{\io}{\iota_{\rm g}}
\newcommand{\smallqed}{{\tiny ($\Box$)}}

\usepackage[round]{natbib}
\usepackage{enumerate}

\author[Bujt{\'a}s et al.]{Csilla Bujt{\'a}s\affiliationmark{1,2}\thanks{Supported by the Slovenian Research and Innovation Agency (ARIS) under the grants P1-0297, N1-0355, J1-70045.}
  \and Tanja Dravec\affiliationmark{3,2}\thanks{Supported by the Slovenian Research and Innovation Agency (ARIS) under the grants P1-0297, N1-0431.}
    \and Michael A. Henning\affiliationmark{4}
    \and Sandi Klav{\v{z}}ar\affiliationmark{1,2,3}\thanks{Supported by the Slovenian Research and Innovation Agency (ARIS) under the grants P1-0297, N1-0285, N1-0355, N1-0431, J1-70045.}
  }
\title{Bounds on the game isolation number and exact values for paths and cycles}

\affiliation{
  Faculty of Mathematics and Physics, University of Ljubljana, Ljubljana, Slovenia\\
  Institute of Mathematics, Physics and Mechanics, Ljubljana, Slovenia\\
  Faculty of Natural Sciences and Mathematics, University of Maribor, Maribor, Slovenia\\
  Department of Mathematics and Applied Mathematics, University of Johannesburg, Johannesburg, South Africa}
\begin{document}
\publicationdata{vol. 28:2}{2026}{33}{10.46298/dmtcs.16132}{2025-07-29; 2025-07-29; 2026-04-02}{2026-05-25}

\maketitle

\begin{abstract}
  The isolation game is played on a graph $G$ by two players  who take turns playing a vertex such that if $X$ is the set of already played vertices, then a vertex can be selected only if it dominates a vertex from a nontrivial component of  $G \setminus N_G[X]$, where $N_G[X]$ is the set of vertices in $X$ or adjacent to a vertex in $X$. Dominator wishes to finish the game with the minimum number of played vertices, while Staller has the opposite goal. The game isolation number $\iota_{\rm g}(G)$ is the number of moves in the Dominator-start game where both players play optimally. If Staller starts the game the invariant is denoted by $\iota_{\rm g}'(G)$. In this paper, $\iota_{\rm g}(C_n)$, $\iota_{\rm g}(P_n)$, $\iota_{\rm g}'(C_n)$, and $\iota_{\rm g}'(P_n)$ are determined for all $n$. It is proved that there are only two graphs that attain equality in the upper bound $\iota_{\rm g}(G) \le \frac{1}{2}|V(G)|$, and that there are precisely eleven graphs which attain equality in the upper bound $\iota_{\rm g}'(G) \le \frac{1}{2}|V(G)|$. For trees $T$ of order at least three it is proved that $\iota_{\rm g}(T) \le \frac{5}{11}|V(T)|$. A new infinite family of graphs $G$ is also constructed for which $\iota_{\rm g}(G) = \iota_{\rm g}'(G) = \frac{3}{7}|V(G)|$ holds. 
\end{abstract}

\section{Introduction}
\label{sec:intro}

Let $G=(V(G),E(G))$ be a graph and let $X\subseteq V(G)$. The closed neighborhood of $X$ is defined as $N_G[X] = \cup_{x\in X}N_G[x]$, where $N_G[x]$ is the closed neighborhood of $x$ in $G$. We say that $X$ 
is an {\em isolating set} of $G$ if $V(G)\setminus N_G[X]$ is an independent set in $G$. The minimum cardinality among all isolating sets in $G$ is the \emph{isolation number} $\iota(G)$ of $G$. The concept was introduced (in a more general setting) in 2017 by~\cite{ch-2017}. For an insight into current developments in the field, we recommend the following articles~\cite{borg-2025, boyer-2024, boyer-2025a, boyer-2025b, chen-2025, lemanska-2024, zhang-2024a, zhang-2024b} and references therein.  

The domination game was introduced by~\cite{bkr-2010}. It was later extensively investigated, the summary of these researches up to 2021  is collected in the book of~\cite{bhkr-2021}. Recent developments on the domination game include the papers by~\cite{bujtas-2022, james-2023, portier-2024, versteegen-2024}. 

Moreover, several variations of the domination game were also investigated, notably the total domination game by~\cite{hkr-2015}, the fractional domination game by~\cite{bt-2019}, the connected domination game by~\cite{bfs-2019}, two variants of the paired-domination game by~\cite{gh-2023, hh-2019}, the competition-independence game of~\cite{ph-sl-1}, and the indicated domination game of~\cite{bresar-2024}. 

More recently, \cite{bdj-2024} introduced a game counterpart of the isolation number following the ideas of the classical domination game (\cite{bkr-2010}). The \emph{isolation game} is played on a graph $G$ by Dominator and Staller, who take turns selecting/playing a vertex from $G$ while obeying the rule that if $X$ is the set of already played vertices, then a vertex can be selected only if it dominates a vertex that belongs to a nontrivial component of the graph $G - N_G[X]$. Such a vertex $x$ is {\em playable}, and a move selecting $x$ is {\em legal}. The game ends when no playable vertex exists. When the game ends, the set of vertices selected forms an isolating set of $G$. As in the domination game, 
Dominator wishes to finish the game by playing the fewest number of vertices, while Staller wishes to delay the process as much as possible. Thus, Dominator seeks to minimize the size of the chosen set while Staller tries to make it as large as possible. If Dominator starts the game, we speak of a {\em D-game}, otherwise it is an {\em S-game}. The {\em game isolation number}, $\io(G)$, and the {\em Staller-start game isolation number}, $\io'(G)$, are the number of vertices selected in the D-game and the S-game, respectively, provided that both players play optimally. 

As mentioned above, the isolation number of a graph was introduced in a more general frame, the same holds for the isolation game defined by~\cite{bdj-2024}. A few results, such as the Continuation Principle and the fact that for each graph $G$ we have $|\io(G) - \io'(G)|\le 1$, have been proved by~\cite{bdj-2024} for the general situation, while a special focus has been given on the isolation game. It was proved that $\io(G)$ is at most half the order for any graph $G$ and conjectured that the actual sharp upper bound is the ceiling of three-sevenths of the order. The isolation number of paths $P_n$ was also determined for the cases when $n \equiv i \bmod 5$, $i\in \{1, 2, 3\}$. 

In Section~\ref{sec:paths-and-cycles}, we determine the exact values for the (Staller-start) game isolation number for paths and cycles. Section~\ref{sec:1/2} focuses on graphs whose game isolation number equals one half of their order. We first show that $K_2$ and $C_6$ are the only graphs that attain the bound $\io(G) \leq \frac{1}{2}|V(G)|$. We then prove that $\io'(G) \leq \frac{1}{2}|V(G)|$ holds for every graph $G$ and identify all connected graphs for which equality holds. There are exactly eleven such graphs, all with order at most 10. In Section~\ref{sec:tree-bound}, we show that the $1/2$-upper bound on the game isolation number in general graphs can be improved to $5/11$ for trees. Regarding the $3/7$-conjecture posed by~\cite{bdj-2024}, we present in Section~\ref{sec:3/7} a new infinite family of graphs with minimum degree at least 2 that achieve this bound. 

\section{Game isolation number of paths and cycles}
\label{sec:paths-and-cycles}

\cite{bdj-2024} established the following results on the game isolation number of paths.

\begin{theorem}{\emph{\cite[Corollary 4.5]{bdj-2024}}} 
\label{thm:DRpaths-1}
If $n \ge 6$, then
\begin{displaymath}
    \left\lceil\frac{2n}{5}\right\rceil-1\le \io(P_n)\le\io'(P_n) \le  \left\lfloor\frac{2n+2}{5}\right\rfloor.
\end{displaymath}
\end{theorem}

\begin{theorem}{\emph{\cite[Corollary 4.6]{bdj-2024}}}
\label{thm:DRpaths-2}
Let $n \ge 6$.  If $n \bmod 5 \in \{1, 2, 3\}$, then
\begin{displaymath}
    \io(P_n)=\io'(P_n)= \left\lfloor\frac{2n+2}{5}\right\rfloor\,.
\end{displaymath}
\end{theorem}

In the first main result of this section we determine the game isolation number for all cycles. After that, in our second main result we use some tools from the proof for cycles to round Theorems~\ref{thm:DRpaths-1} and~\ref{thm:DRpaths-2} by determining the exact values for all paths.

The result for cycles reads as follows. 

\begin{theorem} \label{thm:cycle}
For every $n \ge 4$,
  \begin{displaymath}
      \io(C_n) =  \left \{
    \begin{array}{ll}
	2 \left\lceil \frac{n}{5}\right\rceil; &\text{if~} n \equiv 0 \bmod 5,\\\\
	2 \left\lceil \frac{n}{5}\right\rceil-1; &\text{otherwise}. \\
    \end{array}
    \right.
  \end{displaymath}
    \medskip
\begin{displaymath}
    \io'(C_n) =  \left \{
    \begin{array}{ll}
    2 \left\lfloor \frac{n}{5}\right\rfloor+1; &\text{if~} n \equiv 4 \bmod 5, \\\\
	2 \left\lfloor \frac{n}{5}\right\rfloor; &\text{otherwise}.\\
	
    \end{array}
    \right.
\end{displaymath}
    \end{theorem}

\begin{proof}
For $C_4$ and $C_5$, the formulas can be checked directly. In the continuation, we may assume $n \ge 6$ and prove lower and upper bounds on the invariants by observing strategies for Staller and Dominator, respectively in the subsequent two claims. 
\begin{claim} \label{claim:A}
 For every $n \ge 6$,
     \begin{equation} \label{eq:1}
     \min\left\{2 \left\lceil \frac{n}{5}\right\rceil, 2 \left\lceil \frac{n-4}{5}\right\rceil +1\right\} \le \io(C_n)
     \end{equation}
     and
     \begin{equation} \label{eq:2}
     \min\left\{2 \left\lceil \frac{n-3}{5}\right\rceil, 2 \left\lceil \frac{n-4}{5}\right\rceil +1\right\} \le \io'(C_n).
     \end{equation}    
\end{claim}
\proof We consider the following strategy of Staller. In each of her moves, she plays a (playable) vertex that is a neighbor of an already played vertex. Observe that it is always possible when the game is not over. By a run $B_i$ we mean a maximal sequence of played vertices along the cycle.  Let $B_1, \dots, B_\ell$ be the runs obtained when the game is over. According to Staller's strategy, she never increases the number of runs except in the first move in the S-game. Hence, the number of Dominator's moves is an upper bound on $\ell$ and $\ell - 1$ in the D-game and S-game, respectively.
 Between any two consecutive runs there are at most three unplayed vertices. Let  \begin{math} t = \sum_{i=1}^{\ell} |B_i|,  \end{math} which is exactly the number of played vertices during the game.
 \medskip

We first prove (\ref{eq:1}) and assume that a D-game is played on $C_n$. If $t$ is even, Staller finishes the game and the number of runs can be estimated as $\ell \le t/2 $. Therefore,
 \begin{displaymath}
      n \le t + 3 \,\frac{t}{2}= 5\, \frac{t}{2}
 \end{displaymath}
 and, as $\frac{t}{2}$ is an integer, we have \begin{math}
     \frac{t}{2} \ge \left\lceil \frac{n}{5}\right\rceil.
 \end{math}  Hence  
 we may conclude \begin{math}
     2 \left\lceil \frac{n}{5} \right\rceil \le t
 \end{math}  when $t$ is even. If $t$ is odd, then Dominator plays $\frac{t+1}{2}$ vertices and \begin{math}
     \ell \le \frac{t+1}{2}.
 \end{math}  This implies
 \begin{displaymath}
     n \le t + 3 \,\frac{t+1}{2}= 5\,\frac{t-1}{2}+4.
 \end{displaymath}
  Since $\frac{t-1}{2}$ is an integer, we get
\begin{math}
    2 \left\lceil \frac{n-4}{5} \right\rceil +1 \le t
\end{math} 
  for the case of $t$ is odd. Putting together the cases with odd and even $t$, we obtain inequality (\ref{eq:1}).
\medskip

For the $S$-game, we can prove (\ref{eq:2}) in an analogous way. Here, the first move of Staller creates a new run, but later she does not increase the number of runs. Hence, the number of Dominator's moves is an upper bound on $\ell-1$. If $t$, the number of played vertices, is an even number, then \begin{math}
    \ell \le \frac{t}{2}+1
\end{math}  and
\begin{displaymath}
    n \le t + 3 \left(\frac{t}{2}+1\right)= 5\, \frac{t}{2} +3.
 \end{displaymath}
 Therefore, \begin{math}
     2 \left\lceil \frac{n-3}{5} \right\rceil \le t
 \end{math} holds if $t$ is even. If $t$ is odd, then Staller finishes the game, and we get
\begin{displaymath}
    n \le t + 3 \left(\frac{t-1}{2}+1\right)= 5\,\frac{t-1}{2}+4
\end{displaymath} 
 and conclude
 \begin{math}
     2 \left\lceil \frac{n-4}{5} \right\rceil +1 \le t.
 \end{math} From the odd and even cases, we get (\ref{eq:2}) as stated. \smallqed

\begin{claim} \label{claim:B}
  For every $n \ge 6$, the following inequalities hold.
  \begin{enumerate}[(a)]
      \item
      If \begin{math} n \not\equiv 0 \bmod{5} \end{math}, then \begin{math}\io(C_n) \le 2 \left\lceil \frac{n}{5}\right\rceil-1 \end{math}.
      If \begin{math} n \equiv 0 \bmod{5} \end{math}, then \begin{math} \io(C_n) \le 2 \left\lceil \frac{n}{5}\right\rceil \end{math}.
      \item If \begin{math} n \not\equiv 4 \bmod{5} \end{math}, then \begin{math} \io'(C_n) \le 2 \left\lfloor \frac{n}{5}\right\rfloor \end{math}. If \begin{math} n \equiv 4 \bmod{5} \end{math}, then \begin{math} \io'(C_n) \le 2 \left\lfloor \frac{n}{5}\right\rfloor+1 \end{math}.
  \end{enumerate}     
\end{claim}
\proof To prove the upper bounds stated in the claim, we consider the following strategy of Dominator. If it is not the first move in the game, then Dominator chooses a vertex $v_i$ such that $v_{i-4}$ is a played vertex and $v_{i-3}, v_{i-2}, v_{i-1}$ are all unplayed. If the game is not over, then he is able to follow this strategy. We further note that if there is a playable vertex $v_j$ in $C_n$ (at any moment in the game), then there is a sequence of at least four consecutive playable vertices that includes $v_j$. Consequently, if $C_n$ contains at most seven playable vertices before a move of Dominator, then these vertices form a sequence, and Dominator can choose one of them to finish the game with this move. If the number of playable vertices is eight, then Dominator can follow his strategy, which ensures that the next move played by Staller will finish the game.  For both the D-game and S-game, we count the playable and unplayable vertices after each move in $C_n$ to prove the upper bounds. Let $u_i$ denote the number of vertices made unplayable by the $i^{\rm th}$ move in the game.
  \medskip

  (a) Assume that a D-game is played. The first move of Dominator makes exactly one vertex (the played vertex itself) unplayable, $u_1=1$. Then, for Staller's move, we have $u_2 \ge 1$. As Dominator's strategy ensures, his every later move makes at least four vertices unplayable and hence $u_{2j-1}\ge 4$ for every $j \ge 2$. For Staller's moves, we have $u_{2j} \ge 1$. In general, there are 
  \begin{displaymath}
      U_{2i}=\sum_{s=1}^{2i} u_s \ge 2 + 5(i-1)=5i-3
  \end{displaymath}
    unplayable vertices after $2i$ moves in the game.
 \begin{itemize}
     \item If $n=5k$, then \begin{math}
          U_{2k-2} \ge 5k-8
     \end{math} and at most eight playable vertices remain before Dominator's move. He can make sure that the game finishes with the $(2k)^{\rm th}$ move or sooner. Therefore, \begin{math}
         \io(C_n) \le 2k = 2 \left\lceil \frac{n}{5}\right\rceil
     \end{math}  as stated.
     \item If $n=5k+r$ with $r \in \{1,\dots, 4\}$, then $U_{2k} \ge 5k-3$ implies that at most $r+3 \le 7$ vertices are playable before Dominator's move and he can end the game with the $(2k+1)^{\rm st}$ move (if it was not finished earlier). Therefore, \begin{math}
         \io(C_n) \le 2k+1 = 2 \left\lceil \frac{n}{5}\right\rceil-1. 
     \end{math}
 \end{itemize}

(b) In an S-game on $C_n$ Dominator's strategy implies $u_1=1$, $u_2 \ge 4$, and in general, $u_{2j-1} \ge 1$ and $u_{2j} \ge 4$. After $2i-1$ moves, the number of unplayable vertices, $U_{2i-1}$, is at least $5i-4$.
\begin{itemize}
    \item If $n=5k+4$, then $U_{2k-1} \ge 5k-4$. Thus, at most eight playable vertices remain and the game finishes with the $(2k+1)^{\rm st}$ move or sooner. We conclude that  \begin{math} \io'(C_n) \le 2k+1 = 2 \left\lfloor \frac{n}{5}\right\rfloor+1 \end{math}.
    \item If $n=5k+r$ with $r \in \{0,\dots, 3\}$, then $U_{2k-1} \ge 5k-4$ implies that the number of playable vertices is at most $4+r \le 7$. It follows that Dominator can finish the game with his next move and \begin{math}
        \io'(C_n) \le 2k=  2 \left\lfloor \frac{n}{5}\right\rfloor.
    \end{math}  \smallqed
\end{itemize}
\medskip

To complete the proof of the theorem, we compare the lower and upper bounds in Claims~\ref{claim:A} and \ref{claim:B}.

If $n=5k$, then
\begin{displaymath}
    2 \left\lceil \frac{n}{5}\right\rceil =2k < 2k+1 =
2 \left\lceil \frac{n-4}{5}\right\rceil +1,
\end{displaymath} 
and (\ref{eq:1}) gives \begin{math}
    2 \left\lceil \frac{n}{5}\right\rceil \le \io(C_n).
\end{math}  Then, the upper bound in Claim~\ref{claim:B}(a) gives the desired equality. If $n=5k+r$ with $r \in \{1, \dots, 4\}$, then
\begin{displaymath}
    2 \left\lceil \frac{n}{5}\right\rceil =2k+2 > 2k+1 =
2 \left\lceil \frac{n-4}{5}\right\rceil +1,
\end{displaymath} 
and (\ref{eq:1}) implies $2k+1 \le \io(C_n)$. On the other hand, Claim~\ref{claim:B}(a) establishes \begin{math}
    \io(C_n) \le 2 \left\lceil \frac{n}{5}\right\rceil-1= 2k+1.
\end{math}  We may therefore infer \begin{math}
    \io(C_n) =  2 \left\lceil \frac{n}{5}\right\rceil-1.
\end{math} 

For $\io'(C_n)$, we first consider the case of $n=5k+4$ and observe
\begin{displaymath}
    2 \left\lceil \frac{n-3}{5}\right\rceil =2k+2 > 2k+1 =
2 \left\lceil \frac{n-4}{5}\right\rceil +1 .
\end{displaymath} 
  Inequality (\ref{eq:2}) then yields $2k+1 \le \io'(C_n)$ which matches the upper bound in Claim~\ref{claim:B}(b) and proves \begin{math}
      \io'(C_n) = 2 \left\lfloor \frac{n}{5}\right\rfloor+1. 
  \end{math}
 If $n=5k +r$ with $r \in \{0,\dots, 3\}$, then
 \begin{displaymath}
    2 \left\lceil \frac{n-3}{5}\right\rceil =2k < 2k+1 = 2 \left\lceil \frac{n-4}{5}\right\rceil +1 
 \end{displaymath} 
and (\ref{eq:2}) gives $2k \le \io'(C_n)$. By Claim~\ref{claim:B}(b) we know that \begin{math}
  \io'(C_n) \le 2k = 2 \left\lfloor \frac{n}{5}\right\rfloor.  
\end{math} 
This finishes the proof of the theorem.
\end{proof}

We now complete Theorems~\ref{thm:DRpaths-1} and \ref{thm:DRpaths-2} by establishing the exact values of $\io(P_n)$ and $\io'(P_n)$ when $n \equiv 0 \bmod 5$, and  $n \equiv 4 \bmod 5$.   

\begin{theorem}
\label{thm:paths}
Let $n \ge 6$.  If $n \equiv 0 \bmod 5$, or $n \equiv 4 \bmod 5$, then
\begin{displaymath}
    \io(P_n)= \left\lceil\frac{2n}{5}\right\rceil -1 \qquad {\rm and } \qquad \io'(P_n)= \left\lfloor\frac{2n+2}{5}\right\rfloor.
\end{displaymath}
\end{theorem}

\begin{proof}    
 We first prove that the lower bound 
\begin{math}
     \left\lceil\frac{2n}{5}\right\rceil -1 \leq \io(P_n)
\end{math} in Theorem~\ref{thm:DRpaths-1} is tight when $n=5k$ or $n=5k+4$. Consider a D-game on $P_n \colon v_1\dots v_n$ when Dominator applies the strategy described in the proof of Theorem~\ref{thm:cycle},  Claim~\ref{claim:B}. More exactly, Dominator plays $v_3$ as his first move, and later he chooses to play a vertex $v_i$ such that $v_{i-4}$ is played and $v_{i-3}, v_{i-2}, v_{i-1}, v_{i} $ are all unplayed. If there is no such vertex $v_i$, then either the game is over or the last three vertices of the path are unplayed, and the only edge which is not isolated is $v_{n-1}v_n$. Let $t$ denote the number of moves in this game.

Let $u_j$ denote the number of vertices made unplayable by the $j^{\rm th}$ move of the game and let \begin{math}
    U_j= \sum_{s=1}^{j} u_s 
\end{math} be the total number of unplayable vertices after the $j^{\rm th}$ move. By Dominator's strategy, $u_1 \ge 3$, and $u_{2j-1} \ge 4$ for every $j \ge 2$ if $t > 2j-1$. Therefore, \begin{math}
     U_{2j} \ge 5j-1 
\end{math} holds if $t > 2j$. 

If $n=5k$, then \begin{math}
    U_{2k-2}\ge 5k-6.
\end{math}  Hence, after $2k-2$ moves in the game, at most six vertices are playable. Then these six vertices must be consecutive on the path, and Dominator can finish the game with his next move. It gives \begin{math}
    \io(P_n) \leq t \leq 2k-1 = \left\lceil\frac{2n}{5}\right\rceil -1
\end{math}
and proves the formula for $n=5k$. If $n=5k+4$, then we can give a similar reasoning. In this case, $U_{2k}\ge 5k-1 $ that is, at most five vertices are playable after $2k$ moves. These are consecutive vertices, and Dominator can finish the game with his next move. We may conclude
\begin{math}
    \io(P_n) \leq t \leq 2k+1 = \left\lceil\frac{2n}{5}\right\rceil -1
\end{math} 
and equality follows by Theorem~\ref{thm:DRpaths-1}.
\medskip

Consider now an S-game on $P_n$. To prove a lower bound on $\io'(P_n)$, we suppose that Staller plays according to the strategy described in 
the proof of Theorem~\ref{thm:cycle}, Claim~\ref{claim:A}. Her first move is $v_1$, and later she always plays a vertex next to an already played vertex. By this strategy, if $\ell$ is the number of runs (maximal sequences of consecutive played vertices) at the end of the game, then Dominator played at least $\ell-1$ vertices. Recall also that there are at most three unplayed vertices between two consecutive runs. By Staller's strategy, $v_1$ belongs to a run. At the other end of the path, at most two unplayed vertices may remain. Again, let $t$ denote the number of moves in the game.

If $t$ is even, then Dominator played $t/2$ vertices and
\begin{displaymath}
    n \leq t+ 3 \ell-1 \leq t+ 3\left( \frac{t}{2}+1\right) -1 = 5\, \frac{t}{2}+2.
\end{displaymath} 
Since $\frac{t}{2}$ is an integer, it is equivalent to \begin{math}
     2\left\lceil\frac{n-2}{5}\right\rceil \leq t.
\end{math} 
If $t$ is odd, then Dominator played $(t-1)/2 $ vertices and 
\begin{displaymath}
    n \leq t+ 3 \ell-1 \leq t+ 3\left( \frac{t-1}{2}+1\right) -1 = 5\, \frac{t-1}{2}+3
\end{displaymath} 
follows. Hence we get \begin{math}
    2\left\lceil\frac{n-3}{5}\right\rceil +1\leq t. 
\end{math} From the two cases, we derive
\begin{equation} \label{eq:3}
  m= \min\left\{ 2\left\lceil\frac{n-2}{5}\right\rceil, 2\left\lceil\frac{n-3}{5}\right\rceil +1\right\} \leq t \leq \io'(P_n).
 \end{equation}
 
Finally, we compare~(\ref{eq:3}) with the upper bound \begin{math}
    \io'(P_n) \leq \left\lfloor\frac{2n+2}{5}\right\rfloor
\end{math}  from Theorem~\ref{thm:DRpaths-1}.
If $n=5k$, then $m=2k$, and it equals the upper bound. If $n=5k+4$, then  $m=2k+2$, which matches the upper bound again. We conclude that \begin{math}
   \io'(P_n) = \left\lfloor\frac{2n+2}{5}\right\rfloor 
\end{math} holds in both cases as stated. 
\end{proof}

Combining Theorems~\ref{thm:paths} and~\ref{thm:DRpaths-2} and verifying the formulas directly for $n \in \{1,\dots,5\}$, the complete result for paths can be stated as follows. 

\begin{theorem} \label{thm:path-new}
For every $n \ge 1$, we have 
\begin{displaymath}
    \io(P_n) =  \left \{
    \begin{array}{ll}
\left\lfloor\frac{2n+1}{5}\right\rfloor-1; &\text{if~} n \equiv 0 \bmod{5},\\\\	
       \left\lfloor\frac{2n+1}{5}\right\rfloor; &\text{otherwise},
	 \\
    \end{array}
    \right.
\end{displaymath}
   and
\begin{displaymath}
     \io'(P_n) =  \left\lfloor\frac{2n+2}{5}\right\rfloor. 
\end{displaymath}
   \end{theorem}

\section{Graphs that attain the 1/2-bound}
\label{sec:1/2}

It was proved by~\cite{bdj-2024} that the game isolation number of a graph is at most one-half its order. We characterize the graphs achieving equality in this upper bound on the game isolation number. For this purpose, we follow the proof technique employed in~\cite{bdj-2024}, with the modification that we adopt a greedy strategy for Dominator. 

For completeness, we include the key strategy of the proof in~\cite{bdj-2024}. During the course of the game, when $S$ is the set of selected vertices, a vertex $v$ is considered \emph{marked} if $v \in N_G[S]$ or if all neighbors of $v$ are in $N_G[S]$. Note that a marked vertex may or may not be playable and a vertex is not playable if and only if its entire closed neighborhood is marked. 

\begin{theorem}
\label{thm:1/2bound}
If $G$ is a connected graph of order~$n$, then $\io(G) \le \frac{1}{2}n$, with equality if and only if $G \in \{K_2,C_6\}$.  
\end{theorem}
\begin{proof}
 Let $G$ be a connected graph of order~$n$. If $n = 1$, then $\io(G) = 0$, while if $n = 2$, then $G = K_2$ and $\io(G) = 1 = \frac{1}{2}n$. Hence we may assume that $n \ge 3$, for otherwise the desired result is immediate. Following the proof strategy employed in~\cite{bdj-2024}, we weight the vertices in the graph $G$ as follows. At the beginning all vertices are unmarked and each vertex has weight~$1$. During the game, when a vertex is marked its weight drops to~$0$. The {\em weight} $\w(G)$ of $G$ is the sum of weights of all vertices of $G$. Thus at the beginning of the game $\w(G)=n$ and when the game ends the weight is zero. Let $M$ be the set of marked vertices in some stage of the game.

We adopt a refinement of Dominator's strategy as given in~\cite{bdj-2024} in that Dominator now adopts a greedy strategy and plays a vertex that results in a maximum weight decrease. At the start of the game, Dominator therefore plays a vertex of maximum degree~$\Delta(G)$, thereby marking at least $\Delta(G) + 1$ vertices and decreasing the weight by at least $\Delta(G) + 1$ on his first move. On each of his subsequent moves, if there exists a component $C$ of $G-M$ of order at least~$3$, then Dominator can always decrease the weight by at least~$3$ by playing a vertex of degree at least~$2$ in $C$. Every move of Staller decreases the weight of $G$ by at least~$1$ since she must dominate at least one (unmarked) vertex in $G - M$.

In the final stage of the game when every component in $G-M$ is a $K_2$-component, every move decreases the weight by at least~$2$. By playing a vertex in a $K_2$-component of $G - M$, we note that Staller can play a move that decreases the weight by exactly~$2$. However according to Dominator's greedy strategy, he plays a marked vertex that is adjacent to the maximum number of $K_2$-components in $G - M$, thereby decreasing the weight by at least~$2$. 

In summary, Dominator's greedy strategy implies all his moves decrease the weight by at least~$3$, except possibly in the final stage of the game when every component in $G-M$ is a $K_2$-component, in which case all his subsequent moves decrease the weight by at least~$2$ and all subsequent moves of Staller decrease the weight by exactly~$2$. Moreover if at least one component of $G-M$ is of order at least~$3$, then Staller's move decreases the weight of $G$ by at least~$1$. Hence every move of Dominator decreases the weight by at least~$2$, while each of Dominator's move together with Staller's reply to his move decrease the weight by at least~$4$. Therefore on average every played move decreases the weight by at least~$2$, implying that the game finishes after at most $\frac{1}{2}n$ moves.

By our earlier assumption, recall that $n \ge 3$, and so $\Delta(G) \ge 2$. Suppose that $\Delta(G) \ge 3$. In this case, Dominator's first move decreases the weight by at least~$4$. Since Staller's move decreases the weight by at least~$1$, after the first two moves are played at least five vertices are marked. Thereafter, Dominator's greedy strategy guarantees that at most~$\frac{1}{2}(n-5)$ additional moves are needed to finish the game. Hence in this case the  game finishes after at most \begin{math}
    2 + \frac{1}{2}(n-5) = \frac{1}{2}(n-1) 
\end{math} moves, and so $\io(G) < \frac{1}{2}n$.

Hence we may assume that $\Delta(G) = 2$. Thus, $G$ is a path $P_n$ or a cycle $C_n$. Suppose that $G = P_n$ (and $n \ge 3$). Then by Theorem~\ref{thm:path-new}, we have \begin{math}
    \io(G) \le \frac{2}{5}(n+1)  < \frac{1}{2}n.
\end{math}  Hence if $G = P_n$ and $n \ge 3$, then $\io(G) < \frac{1}{2}n$.

Suppose finally that $G = C_n$. If $n \in \{3,4\}$, then $\io(G) = 1 < \frac{1}{2}n$. If $n = 5$, then $\io(G) = 2 < \frac{1}{2}n$. If $n = 6$, then $\io(G) = 3 = \frac{1}{2}n$. Hence we may assume that $n \ge 7$. By Theorem~\ref{thm:cycle}, we now infer the following. If $n \equiv 0 \bmod{5}$, then \begin{math}
    \io(C_n) =  \frac{2}{5}n < \frac{1}{2}n.
\end{math}  If $n \equiv 1 \bmod{5}$, then \begin{math}
    \io(C_n) =  \frac{2}{5}(n+4) - 1 < \frac{1}{2}n. 
\end{math} If $n \equiv 2 \bmod{5}$, then \begin{math}
    \io(C_n) =  \frac{2}{5}(n+3) - 1 < \frac{1}{2}n. 
\end{math} If $n \equiv 3 \bmod{5}$, then \begin{math}
    \io(C_n) =  \frac{2}{5}(n+2) - 1 < \frac{1}{2}n. 
\end{math} If $n \equiv 4 \bmod{5}$, then \begin{math}
    \io(C_n) =  \frac{2}{5}(n+1) - 1 < \frac{1}{2}n. 
\end{math} In all cases, when $n \ge 7$ we have $\io(C_n) < \frac{1}{2}n$. Hence we have shown that $G = K_2$ and $G = C_6$ are the only graphs that satisfy $\io(C_n) = \frac{1}{2}n$.   
\end{proof}

\begin{figure}[htb]
\begin{center}
\begin{tikzpicture}[scale=.85,style=thick,x=0.75cm,y=0.75cm]
\def\vr{3.0pt} 
%
\path (0,1) coordinate (a1);
\path (1,1) coordinate (a2);
\path (3,1) coordinate (b1);
\path (4,1) coordinate (b2);
\path (5,1) coordinate (b3);
\path (6,1) coordinate (b4);
\path (8,0) coordinate (c1);
\path (9,1) coordinate (c2);
\path (9,2) coordinate (c3);
\path (10,0) coordinate (c4);
\path (12,1) coordinate (d1);
\path (13,1) coordinate (d2);
\path (14,1) coordinate (d3);
\path (15,1) coordinate (d4);
\path (16,1) coordinate (d5);
\path (17,1) coordinate (d6);
\path (18,1) coordinate (d7);
\path (15,2) coordinate (d);
\draw (a1)--(a2);
\draw (b1)--(b2)--(b3)--(b4);
\draw (c3)--(c2)--(c1)--(c4)--(c2);
\draw (d1)--(d2)--(d3)--(d4)--(d5)--(d6)--(d7)--(d1);
\draw (d)--(d4);
\draw (a1) [fill=white] circle (\vr);
\draw (a2) [fill=white] circle (\vr);
\draw (b1) [fill=white] circle (\vr);
\draw (b2) [fill=white] circle (\vr);
\draw (b3) [fill=white] circle (\vr);
\draw (b4) [fill=white] circle (\vr);
\draw (c1) [fill=white] circle (\vr);
\draw (c2) [fill=white] circle (\vr);
\draw (c3) [fill=white] circle (\vr);
\draw (c4) [fill=white] circle (\vr);
\draw (d1) [fill=white] circle (\vr);
\draw (d2) [fill=white] circle (\vr);
\draw (d3) [fill=white] circle (\vr);
\draw (d4) [fill=white] circle (\vr);
\draw (d5) [fill=white] circle (\vr);
\draw (d6) [fill=white] circle (\vr);
\draw (d7) [fill=white] circle (\vr);
\draw (d) [fill=white] circle (\vr);
\draw (0.5,-1) node {{\small (a) $F_1$}};
\draw (4.5,-1) node {{\small (b) $F_2$}};
\draw (9,-1) node {{\small (c) $F_3$}};
\draw (15,-1) node {{\small (d) $F_4$}};
\end{tikzpicture}

\vskip 0.5 cm 
\begin{tikzpicture}[scale=.75,style=thick,x=0.75cm,y=0.75cm]
\def\vr{2.5pt} 
%
\path (0,0) coordinate (e1);
\path (0,2) coordinate (e2);
\path (1,1) coordinate (e3);
\path (2,1) coordinate (e4);
\path (3,1) coordinate (e5);
\path (4,1) coordinate (e6);
\path (5,1) coordinate (e7);
\path (2,2) coordinate (e);
\path (10,0) coordinate (g1);
\path (10,2) coordinate (g2);
\path (11,1) coordinate (g3);
\path (12,1) coordinate (g4);
\path (13,1) coordinate (g5);
\path (14,2) coordinate (g6);
\path (14,0) coordinate (g7);
\path (12,2) coordinate (g);
%

\draw (e3)--(e2)--(e1)--(e3)--(e4)--(e5)--(e6)--(e7);
\draw (e)--(e4);
%
\draw (g3)--(g2)--(g1)--(g3)--(g4)--(g5)--(g6)--(g7)--(g5);
\draw (g)--(g4);

\draw (e1) [fill=white] circle (\vr);
\draw (e2) [fill=white] circle (\vr);
\draw (e3) [fill=white] circle (\vr);
\draw (e4) [fill=white] circle (\vr);
\draw (e5) [fill=white] circle (\vr);
\draw (e6) [fill=white] circle (\vr);
\draw (e7) [fill=white] circle (\vr);
\draw (e) [fill=white] circle (\vr);
\draw (g1) [fill=white] circle (\vr);
\draw (g2) [fill=white] circle (\vr);
\draw (g3) [fill=white] circle (\vr);
\draw (g4) [fill=white] circle (\vr);
\draw (g5) [fill=white] circle (\vr);
\draw (g6) [fill=white] circle (\vr);
\draw (g7) [fill=white] circle (\vr);
\draw (g) [fill=white] circle (\vr);
\draw (2.5,-1) node {{\small (e) $F_5$}};
\draw (12,-1) node {{\small (f) $F_6$}};
\end{tikzpicture}

\vskip 0.5 cm
\begin{tikzpicture}[scale=.75,style=thick,x=0.75cm,y=0.75cm]
\def\vr{2.5pt} 
%
\path (0,0.75) coordinate (h1);
\path (0,1.5) coordinate (h2);
\path (1,0) coordinate (h3);
\path (1,2.25) coordinate (h4);
\path (1,3) coordinate (h5);
\path (1,3.75) coordinate (h6);
\path (2,0.75) coordinate (h7);
\path (2,1.5) coordinate (h8);
\path (4,0.75) coordinate (i1);
\path (4,1.5) coordinate (i2);
\path (5,0) coordinate (i3);
\path (5,2.25) coordinate (i4);
\path (5,3) coordinate (i5);
\path (5,3.75) coordinate (i6);
\path (6,0.75) coordinate (i7);
\path (6,1.5) coordinate (i8);
\path (8,0.75) coordinate (j1);
\path (8,1.5) coordinate (j2);
\path (9,0) coordinate (j3);
\path (9,2.25) coordinate (j4);
\path (9,3) coordinate (j5);
\path (9,3.75) coordinate (j6);
\path (10,0.75) coordinate (j7);
\path (10,1.5) coordinate (j8);
\path (12,0.75) coordinate (k1);
\path (12,1.5) coordinate (k2);
\path (13,0) coordinate (k3);
\path (13,2.25) coordinate (k4);
\path (13,3) coordinate (k5);
\path (13,3.75) coordinate (k6);
\path (14,0.75) coordinate (k7);
\path (14,1.5) coordinate (k8);
\path (16,0.75) coordinate (l1);
\path (16,1.5) coordinate (l2);
\path (17,0) coordinate (l3);
\path (17,2.25) coordinate (l4);
\path (17,3) coordinate (l5);
\path (17,3.75) coordinate (l6);
\path (18,0.75) coordinate (l7);
\path (18,1.5) coordinate (l8);

\draw (h1)--(h2)--(h4)--(h8)--(h7)--(h3)--(h1);
\draw (h4)--(h5)--(h6);
\draw (i1)--(i2)--(i4)--(i8)--(i7)--(i3)--(i1);
\draw (i4)--(i5)--(i6);
\draw (i2) to[out=90,in=180, distance=0.5cm ] (i5);
\draw (j1)--(j2)--(j4)--(j8)--(j7)--(j3)--(j1);
\draw (j5)--(j6);
\draw (j2) to[out=90,in=180, distance=0.5cm ] (j5);
\draw (j8) to[out=90,in=0, distance=0.5cm ] (j5);
\draw (k1)--(k2)--(k4)--(k8)--(k7)--(k3)--(k1);
\draw (k5)--(k6);
\draw (k2) to[out=90,in=180, distance=0.5cm ] (k5);
\draw (k7) to[out=45,in=0, distance=0.85cm ] (k5); 
\draw (l1)--(l2)--(l4)--(l8)--(l7)--(l3)--(l1);
\draw (l4)--(l5)--(l6);
\draw (l2) to[out=90,in=180, distance=0.5cm ] (l5);
\draw (l8) to[out=90,in=0, distance=0.5cm ] (l5);

\draw (h1) [fill=white] circle (\vr);
\draw (h2) [fill=white] circle (\vr);
\draw (h3) [fill=white] circle (\vr);
\draw (h4) [fill=white] circle (\vr);
\draw (h5) [fill=white] circle (\vr);
\draw (h6) [fill=white] circle (\vr);
\draw (h7) [fill=white] circle (\vr);
\draw (h8) [fill=white] circle (\vr);
\draw (i1) [fill=white] circle (\vr);
\draw (i2) [fill=white] circle (\vr);
\draw (i3) [fill=white] circle (\vr);
\draw (i4) [fill=white] circle (\vr);
\draw (i5) [fill=white] circle (\vr);
\draw (i6) [fill=white] circle (\vr);
\draw (i7) [fill=white] circle (\vr);
\draw (i8) [fill=white] circle (\vr);
\draw (j1) [fill=white] circle (\vr);
\draw (j2) [fill=white] circle (\vr);
\draw (j3) [fill=white] circle (\vr);
\draw (j4) [fill=white] circle (\vr);
\draw (j5) [fill=white] circle (\vr);
\draw (j6) [fill=white] circle (\vr);
\draw (j7) [fill=white] circle (\vr);
\draw (j8) [fill=white] circle (\vr);
\draw (k1) [fill=white] circle (\vr);
\draw (k2) [fill=white] circle (\vr);
\draw (k3) [fill=white] circle (\vr);
\draw (k4) [fill=white] circle (\vr);
\draw (k5) [fill=white] circle (\vr);
\draw (k6) [fill=white] circle (\vr);
\draw (k7) [fill=white] circle (\vr);
\draw (k8) [fill=white] circle (\vr);
\draw (l1) [fill=white] circle (\vr);
\draw (l2) [fill=white] circle (\vr);
\draw (l3) [fill=white] circle (\vr);
\draw (l4) [fill=white] circle (\vr);
\draw (l5) [fill=white] circle (\vr);
\draw (l6) [fill=white] circle (\vr);
\draw (l7) [fill=white] circle (\vr);
\draw (l8) [fill=white] circle (\vr);

\draw (1,-1) node {{\small (g) $F_7$}};
\draw (5,-1) node {{\small (h) $F_8$}};
\draw (9,-1) node {{\small (i) $F_{9}$}};
\draw (13,-1) node {{\small (j) $F_{10}$}};
\draw (17,-1) node {{\small (k) $F_{11}$}};
\end{tikzpicture}

\end{center}
\begin{center}
\vskip -0.5 cm
\caption{Graphs in the family~$\cF$}
\label{fig:familyF}
\end{center}
\end{figure}
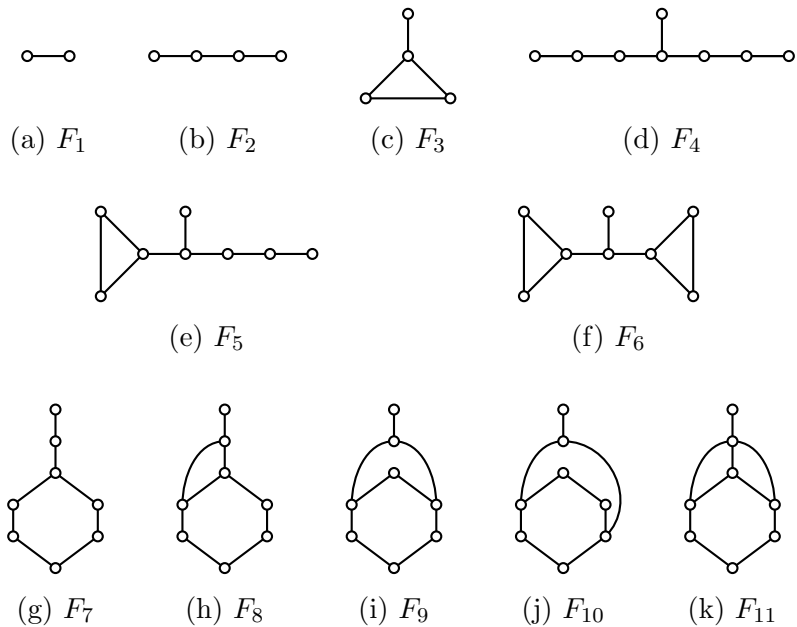

\medskip
We prove next that the Staller-start game isolation number of a graph is at most one-half its order, and we characterize the graphs achieving equality in the $\frac{1}{2}n$-upper bound on the Staller-start game isolation number. For this purpose, let $\cF = \{F_1,F_2,\ldots,F_{11}\}$ be the family of eleven graphs illustrated in Fig.~\ref{fig:familyF}. 

\begin{theorem}
\label{thm:1/2bound-Sgame}
If $G$ is a connected graph of order~$n$, then $\io'(G) \le \frac{1}{2}n$, with equality if and only if $G \in \cF$.
\end{theorem}
  \begin{proof}
Let $G$ be a connected graph of order~$n$. We first prove the desired upper bound. 

\begin{claim}
\label{claim-0}
$\io'(G) \le \frac{1}{2}n$. 
\end{claim}
\proof 
If $n = 1$, then $\io'(G) = 0$, hence assume in the rest that $n \ge 2$. As in the proof of Theorem~\ref{thm:1/2bound}, we show that if Dominator employs his greedy strategy, then the game finishes after at most $\frac{1}{2}n$ moves. Also, let $M$ be the set of marked vertices in some stage of the game. At the start of the game, Staller plays a vertex that marks at least $\delta(G) + 1$ vertices. Therefore her first move decreases the weight by at least~$2$. Thereafter, as shown in the proof of Theorem~\ref{thm:1/2bound}, Staller's subsequent moves decrease the weight by at least~$1$. Moreover in the final stage of the game when every component in $G-M$ is a $K_2$-component, every move of Staller decreases the weight by at least~$2$. As shown in the proof of Theorem~\ref{thm:1/2bound}, Dominator has a strategy that guarantees that all his moves decrease the weight by at least~$3$, except possibly in the final stage of the game when every component in $G-M$ is a $K_2$-component, in which case all his subsequent moves decrease the weight by at least~$2$.

If Dominator's first move decreases the weight by exactly~$2$, then since he follows the greedy strategy, all moves played in the game decrease the weight by exactly~$2$, implying that the game finishes after $\frac{1}{2}n$ moves. Hence we may assume that the first move of Dominator decreases the weight by at least~$3$, for otherwise the desired upper bound follows. Let $i+1$ be the first move of Dominator that decreases the weight by exactly~$2$, where $i \ge 1$. Thus, all moves after Dominator plays his $(i+1)$th move decrease the weight by exactly~$2$. The first move of Staller decreases the weight by at least~$2$, while all her other moves played before Dominator's $(i+1)$th move decrease the weight by at least~$1$. Hence the first $2i+1$ moves played in the game (that is, all moves played immediately before Dominator plays his $(i+1)$th move) decrease the total weight by at least~$2 + i + 3i = 2(2i+1)$. Thus on average every played move decreases the weight by at least~$2$, implying that the game finishes after at most $\frac{1}{2}n$ moves. This establishes the desired upper bound, namely $\io'(G) \le \frac{1}{2}n$.~\smallqed

\medskip 
Suppose next that $\io'(G) = \frac{1}{2}n$. In particular, $n \ge 2$ is even. If $n \in \{2,4\}$, then it is straightforward to verify that in this case $G \in \{F_1,F_2,F_3\}$. Hence we may assume that in the rest of the proof $n \ge 6$. We proceed further with a series of claims establishing some structural properties of the graph $G$. 

\begin{claim}
\label{claim-1}
Staller's first move decreases the weight by exactly~$2$. 
\end{claim}
\proof Suppose that Staller's first move decreases the weight by at least~$3$. As observed earlier, on average every move played after the first move of Staller decreases the weight by at least~$2$, implying that the game finishes after at most \begin{math}
    1 + \frac{1}{2}(n-3) < \frac{1}{2}n
\end{math}  moves, a contradiction.~\smallqed

\medskip 
By Claim~\ref{claim-1}, Staller's first move decreases the weight by exactly~$2$, implying that the first move of Staller is a vertex of degree~$1$. Let $u$ be the vertex played by Staller on her first move, and let $v$ be the (unique) neighbor of $u$ in $G$. Thus at this opening stage of the game, we have $M = \{u,v\}$.

\begin{claim}
\label{claim-2}
Every move of Dominator decreases the weight by~$2$ or~$3$. 
\end{claim}
\proof Suppose that after the first move of Staller, every component of $G - M$ is a $K_2$-component. By assumption $n \ge 6$, and so $G - M$ contains $k \ge 2$ components and $n = 2k + 2$. Dominator now plays the vertex~$v$ as his first move, thereby marking all vertices, and so \begin{math}
    \io'(G) = 2 < k + 1 = \frac{1}{2}n,
\end{math}  a contradiction. Hence, at least one component of $G - M$ has order at least~$3$. The first move of Dominator therefore decreases the weight by at least~$3$.

Suppose that the first move of Dominator decreases the weight by at least~$4$. The second move of Staller decreases the weight by at least~$1$. Hence after the first three moves of the game (that is, after Staller's first and second move and Dominator's first move) at least~$2 + 4 + 1 = 7$ vertices are marked. Thereafter, on average every move played starting with the second move of Dominator decreases the weight by at least~$2$, implying that the game finishes after at most \begin{math}
    3 + \frac{1}{2}(n-7) < \frac{1}{2}n
\end{math}  moves, a contradiction. Hence, the first move of Dominator decreases the weight by exactly~$3$. Analogous argument shows that if Dominator plays a move that decreases the weight by at least~$4$ at any stage in the game, then the game finishes after strictly less than $\frac{1}{2}n$ moves, a contradiction. Hence, every move of Dominator decreases the weight by~$2$ or~$3$.~\smallqed  

\medskip 
Let $U$ denote the set of unmarked vertices, and so \begin{math}
    U = V(G) \setminus M.
\end{math}  

\begin{claim}
\label{claim-3}
At any stage after Staller's first move, the following properties hold. 
\begin{enumerate}[(a)]
\item Every vertex in $M$ has at most three neighbors in~$U$.
\item  Every vertex in $U$ has at most two neighbors in~$U$.
\item  $G - M$ is the disjoint union of paths and cycles.
\item  If it is Dominator's move and $G - M$ contains a path component, then such a component is a $P_2$- or $P_3$-component. 
\end{enumerate}
\end{claim}
\proof If it is Dominator's move and there is vertex in $M$ with four or more neighbors in $U$, then Dominator can play such a vertex and decrease the weight by at least~$4$, contradicting Claim~\ref{claim-2}. Hence if it is Dominator's move, then every vertex in $M$ has at most three neighbors in~$U$.

If it is Dominator's move and there is a vertex in $U$ with three or more neighbors in $U$, then Dominator can play such a vertex and decrease the weight by at least~$4$, contradicting Claim~\ref{claim-2}. Hence if it is Dominator's move, then every vertex in $U$ has at most two neighbors in~$U$.

Immediately after Staller plays her first move, namely the vertex~$u$, we have $M = \{u,v\}$. At this opening stage of the game ($M = \{u,v\}$), properties~(a) and~(b) in the statement of Claim~\ref{claim-3} hold. Since properties~(a) and~(b) are preserved as the game progresses, we infer that properties~(a) and~(b) also hold if it is Staller's move. This proves Parts~(a) and~(b). 

By Part~(b), $G - M$ has maximum degree at most~$2$ and is therefore the disjoint union of paths and cycles. This proves Part~(c).

Suppose that it is Dominator's move and that $G - M$ contains a path component $P \colon v_1 v_2 \ldots v_p$, where $p \ge 4$. In this case, if Dominator plays the vertex~$v_3$, then he marks at least $v_1$, $v_2$, $v_3$, and $v_4$, and therefore decreases the weight by at least~$4$, contradicting Claim~\ref{claim-2}. This proves Part~(d).~\smallqed  

\medskip
We now return to the opening stage of the game, that is, immediately after Staller plays her first move, namely the vertex~$u$, and so $M = \{u,v\}$. Let $G_1, \ldots, G_k$ denote the components of $G - M$. By Claim~\ref{claim-3}(a), the vertex~$v \in M$ is adjacent to at most three vertices in
\begin{math}
    U = V(G) \setminus \{u,v\}, 
\end{math} and so $k \le 3$. By Claim~\ref{claim-3}(c), every component $G_i$ is a path or cycle.

\begin{claim}
\label{claim-4a}
In the initial stage of the game when $M = \{u,v\}$, no component of $G - M$ is a $C_4$-component.
\end{claim}
\proof Suppose that immediately after Staller's first move, there exists a $C_4$-component $C \colon v_1 v_2 v_3 v_4 v_1$ in $G - M$ where $vv_1 \in E(G)$. In this case, if Dominator plays the vertex~$v_1$, then he marks all four vertices in the component $C$. Hence, Dominator can decrease the weight by at least~$4$, contradicting Claim~\ref{claim-2}.~\smallqed 

\begin{claim}
\label{claim-4}
In the initial stage of the game when $M = \{u,v\}$, if $G - M$ contains a path component, then such a component is a $P_3$-component. 
\end{claim}
\proof Suppose that immediately after Staller's first move, there exists a path component $P \colon v_1 v_2 \ldots v_p$ in $G - M$ for some $p \ge 2$. By Claim~\ref{claim-3}(d), we infer that $p \in \{2,3\}$. More generally, every path component in $G - M$ is a $P_2$- or $P_3$-component.

We show that $p = 3$. Suppose, to the contrary, that $p = 2$. Renaming components of $G - M$, we may assume that $G_1 = P$, and so $G_1$ is a $P_2$-component. We may further assume that $vv_1$ is an edge, where recall that $P \colon v_1v_2$. We note that $v$ and $v_2$ may or may not be adjacent. Recall that $k \le 3$, where $k$ denotes the number of components in $G - M$.

Suppose that $k = 3$. In this case, if Dominator plays the vertex~$v$, then he marks both vertices in $G_1$ and he marks at least one vertex in each of $G_2$ and $G_3$, and therefore marks at least four vertices. Hence, Dominator can decrease the weight by at least~$4$, contradicting Claim~\ref{claim-2}. Therefore, $k \le 2$. If $k = 1$, then $n = 4$, contradicting our earlier assumption that $n \ge 6$. Hence, $k = 2$. 

If $G_2$ has odd order, then $G$ would have odd order, contradicting the parity of~$n$ (recall that~$n \ge 6$ is even). Hence, $G_2$ is either a $P_2$-component or a cycle component of even order. If $G_2$ is a $P_2$-component, then $n = 6$. However if Dominator now plays the vertex~$v$, then he marks all four vertices in $G_1$ and $G_2$, implying that the game is finished after Dominator's first move, that is, after strictly less than $\frac{1}{2}n$ moves, a contradiction. Hence, $G_2$ is a cycle of even order. By Claim~\ref{claim-4a}, no component in $G - M$ is a $C_4$-component. Hence, $G_2$ is an even cycle of order at least~$6$. Let $G_2$ be the cycle $w_1 w_2 \ldots w_qw_1$, where $q \ge 6$ is an even integer.

Suppose that $v$ is adjacent to two or more vertices in $G_2$. In this case, as before if Dominator plays the vertex~$v$, then he decreases the weight by at least~$4$, a contradiction. Hence, $v$ is adjacent to exactly one vertex in $G_2$. Renaming vertices if necessary, we may assume that $vw_1$ is an edge.

Suppose that $q = 6$, and so $n = 10$. In this case, Dominator plays as his first vertex the vertex~$v$, thereby marking the two vertices in $G_1$ and the vertex~$w_1$ in $G_2$. Immediately after Dominator's first move, we have \begin{math}
    M = \{u,v,v_1,v_2,w_1\}
\end{math}  and $G - M$ is a $P_5$-component given by $w_2w_3w_4w_5w_6$. If Staller now plays the vertex~$w_i$ on her second move for some $i \in [6]$, then Dominator finishes the game by playing the vertex~$w_{i+3}$ where addition is taken modulo~$6$. Hence, the game is finished in four moves, and so \begin{math}
     \io'(G) \le 4 < \frac{1}{2}n
\end{math} (where recall that $n = 10$). 

Thus, $q \ge 8$ and $q$ is even. In this case, Dominator plays as his first vertex the vertex~$w_1$, thereby marking three new vertices in $G_2$, namely $\{w_1,w_2,w_q\}$. Thus at this stage of the game, we have \begin{math}
    M = \{u,v,w_1,w_2,w_q\}
\end{math}  and $G - M$ consists of two components, namely $G_1$ and the path component $G_2'$ given by $w_3 w_4 \ldots w_{q-1}$. We note that $G_2'$ is a $P_{q-3}$-component where $q \ge 8$. 

Suppose that $q \ge 12$. In this case, whatever move Staller plays as her second move, with the updated set $M$ the graph $G - M$ contains a path component of order at least~$4$, contradicting Claim~\ref{claim-3}(d). Hence, $q \in \{8,10\}$. If $q = 8$, then $n = 12$, but in this case $\io'(G) \le 5 < \frac{1}{2}n$. If $q = 10$, then $n = 14$ and in this case $\io'(G) \le 6 < \frac{1}{2}n$. We therefore deduce that $p = 3$. This completes the proof of Claim~\ref{claim-4}.~\smallqed 

\medskip 
By Claim~\ref{claim-4}, in the initial stage of the game when $M = \{u,v\}$, if $G - M$ contains a path component, then such a component is a $P_3$-component.

\begin{claim}
\label{claim-5}
In the initial stage of the game when $M = \{u,v\}$, if there exist at least two path components of $G - M$, then $G = F_4$. 
\end{claim}
\proof Suppose that immediately after Staller's first move, there exist two or more path components in $G - M$. By Claim~\ref{claim-4}, every path component of $G - M$ is a $P_3$-component. Renaming the components of $G - M$ if necessary, we may assume that $G_1$ and $G_2$ are $P_3$-components. Suppose that $k = 3$. Since $n$ is even, the component $G_3$ has even order, implying that $G_3$ is a cycle of even order. By Claim~\ref{claim-4a}, no component in $G - M$ is a $C_4$-component. Hence, $G_3$ is an even cycle of order at least~$6$. Let $G_3$ be the cycle $w_1 w_2 \ldots w_qw_1$, where $q \ge 6$ is an even integer.

Let $G_1$ be the path $x_1x_2x_3$ and let $G_2$ be the path $y_1y_2y_3$. If $vx_2$ is an edge, then Dominator's first move is~$v$, thereby marking all three vertices in $G_1$ and at least one vertex in each of $G_2$ and $G_3$. Hence in this case, Dominators marks at least five vertices, contradicting Claim~\ref{claim-2}. Hence, $vx_2$ is not an edge. Renaming $x_1$ and $x_3$ if necessary, we may assume that $vx_1 \in E(G)$. Analogously, we may assume that $vy_1 \in E(G)$. Renaming vertices of $G_3$ we may assume that $vw_1 \in E(G)$. 

Now Dominator's first move is~$v$, thereby marking a vertex in each of $G_1$, $G_2$ and $G_3$. In particular, we note that the vertex~$w_1$ is marked. Immediately after Dominator's move, we have $M = \{u,v,x_1,y_1,w_1\}$ and $G - M$ contains three components, namely \begin{math}
    G_1' = G_1 - x_1,
\end{math}  \begin{math}
    G_2' = G_2 - y_1
\end{math}  and \begin{math}
    G_3' = G_3 - w_3.
\end{math}  Moreover, $G_1'$ and $G_2'$ are both $P_2$-components, and $G_3'$ is a $P_{q-1}$-component where $q \ge 6$.

If $q \ge 12$, then whatever move Staller plays as her second move, with the updated set $M$ the graph $G - M$ contains a path component of order at least~$4$, contradicting Claim~\ref{claim-3}(d). Hence, $q \in \{6,8,10\}$. If $q = 6$, then $n = 14$ and we have $\io'(G) \le 6 < \frac{1}{2}n$. If $q = 8$, then $n = 16$ and in this case, $\io'(G) \le 7 < \frac{1}{2}n$. If $q = 10$, then $n = 18$ and in this case,  \begin{math}
    \io'(G) \le 8 < \frac{1}{2}n.
\end{math}  We therefore deduce that $k = 2$. Thus, $G = F_4$.~\smallqed

\medskip
By Claim~\ref{claim-5}, we may assume that in the initial stage of the game when $M = \{u,v\}$, there exists at most one path component in $G - M$.

\begin{claim}
\label{claim-6}
In the initial stage of the game when $M = \{u,v\}$, if $G$ contains a $P_3$-component, then $k = 2$.
\end{claim}
\proof Suppose that immediately after Staller's first move, there exists a $P_3$-component of $G - M$. Since $n \ge 6$ is even, we note that $k \ge 2$. Renaming components if necessary, we may assume that $G_1$ is a $P_3$-component. Suppose, to the contrary, that $k = 3$. By assumption, there exists at most one path component in $G - M$, and so the components $G_2$ and $G_3$ are both cycle components. 

Suppose that $G_2$ or $G_3$ is a $C_3$-component. Then Claim~\ref{claim-2} implies that $v$ has exactly one neighbor in each of the components. Renaming components if necessary, we may assume that $G_2$ is a $C_3$-component. Dominator now plays as his first vertex the vertex~$v$ and proceeding now exactly as in the proof of Claim~\ref{claim-5}, we show that $\io'(G) < \frac{1}{2}n$, a contradiction. Hence, neither $G_2$ nor $G_3$ is a $C_3$-component. By Claim~\ref{claim-4a}, no component in $G - M$ is a $C_4$-component. Hence, each of $G_2$ and $G_3$ is cycle of order at least~$5$.  

As before, Dominator now plays as his first move the vertex~$v$, thereby marking a vertex in each of $G_1$, $G_2$ and $G_3$. Immediately after Dominator's move, we have $M = N_G[v]$ and $G - M$ contains three components, namely one $P_2$-component and two path components that are both of order at least~$4$. Whatever move Staller plays as her second move, with the updated set $M$ the graph $G - M$ contains a path component of order at least~$4$, contradicting Claim~\ref{claim-3}(d). Hence, $k = 2$, as claimed.~\smallqed

\begin{claim}
\label{claim-7}
In the initial stage of the game when $M = \{u,v\}$, if there exists a path component in $G - M$, then $G = F_5$. 
\end{claim}
\proof Suppose that immediately after Staller's first move, there exists a component of $G - M$ that is not a cycle component. Renaming components if necessary, we may assume that $G_1$ is a path component. By Claim~\ref{claim-4}, $G_1$ is a $P_3$-component. Let $G_1$ be the path $v_1v_2v_3$. Since $n \ge 6$ is even, we note that $k \ge 2$.

By assumption, every component of $G - M$ different from $G_1$ is a cycle component. By Claim~\ref{claim-6}, $k = 2$. Since $n$ is even, the component $G_2$ is therefore a cycle component of odd order. Let $G_2$ be the cycle $w_1 w_2 \ldots w_q w_1$, where $q \ge 3$ is odd. Renaming vertices if necessary, we may assume that $vw_1$ is an edge. 

Suppose that $q = 5$, and so $n = 10$. Dominator now plays~$v_2$ as his first move. After Dominator plays the vertex~$v_2$, with the updated set $M$ the graph $G - M$ consists of a $5$-cycle. The game is therefore finished in two more moves, and so \begin{math}
    \io'(G) \le 4 < \frac{1}{2}n,
\end{math}  a contradiction. Hence, $q \ne 5$. 

Suppose that $q = 7$, and so $n = 12$. Dominator now plays as his first vertex the vertex~$v_2$. After Dominator plays the vertex~$v_2$, with the updated set $M$ the graph $G - M$ consists of a $7$-cycle. If Staller plays a vertex of the $7$-cycle on her second move, then Dominator can finish the game in one more move, while if she plays the vertex~$v$ in her second move, then the game finishes in two more moves. Hence, \begin{math}
    \io'(G) \le 5 < \frac{1}{2}n,
\end{math}  a contradiction. Hence, $q \ne 7$.

Suppose that $q \ge 9$ and $q$ is odd. Dominator now plays~$v_2$ as his first move. Whatever move Staller plays as her second move, with the updated set $M$ the graph $G - M$ contains a path component of order at least~$4$, contradicting Claim~\ref{claim-3}(d). 

Hence, $q = 3$, and so $n = 8$. Dominator now plays ~$v$ as his first move. If his move marks all vertices in $G_1$ or $G_2$, then the game is finished after Staller's second move. In this case, \begin{math}
    \io'(G) \le 3 < \frac{1}{2}n,
\end{math}  a contradiction. Hence after Dominator plays the vertex~$v$, with the updated set $M$ the graph $G - M$ contains two $P_2$-components, implying that $G = F_5$.~\smallqed

\medskip
By Claim~\ref{claim-7}, we may assume that in the initial stage of the game when $M = \{u,v\}$, every component in $G - M$ is a cycle component. 

\begin{claim}
\label{claim-8}
In the initial stage of the game when $M = \{u,v\}$, if there exists a $C_3$-component in $G - M$, then $G = F_6$.
\end{claim}
\proof Suppose that immediately after Staller's first move, there exists a $C_3$-component in $G - M$. Since $n \ge 6$ is even, we note that $k \ge 2$. Renaming components if necessary, we may assume that $G_1$ is a $C_3$-component. Let $G_1$ be the cycle $v_1v_2v_3v_1$ where we may assume that $vv_2$ is an edge. If $k = 3$, then proceeding exactly as in the proof of Claim~\ref{claim-6} we obtain a contradiction. Hence, $k = 2$. Since $n$ is even, the component $G_2$ is therefore a $C_q$-component for some $q \ge 3$ where $q$ is odd. If $q \ge 5$, then proceeding exactly as in the proof of Claim~\ref{claim-7} we obtain a contradiction. Hence, $q = 3$. If $v$ is adjacent to two vertices in $G_1$ or in $G_2$, then Dominator plays as his first vertex the vertex~$v$. In this case, with the updated set $M$ the graph $G - M$ consists of one $P_2$-component, and so the game is finished in three moves. Hence, \begin{math}
    \io'(G) \le 3 < \frac{1}{2}n, 
\end{math} a contradiction. Hence, $v$ is adjacent to exactly one vertex in each of $G_1$ and $G_2$. Thus, $G = F_6$.~\smallqed 

\medskip
Suppose that we are in the initial stage of the game when $M = \{u,v\}$. By Claim~\ref{claim-8}, we may assume that no component in $G - M$ is a $C_3$-component. By Claim~\ref{claim-4a}, no component of $G - M$ is a $C_4$-component. Hence, every component in $G - M$ is a cycle component of order at least~$5$. 

\begin{claim}
\label{claim-9}
$k = 1$. 
\end{claim}
\proof Suppose, to the contrary, that $k \ge 2$. Recall that $G_1, \ldots, G_k$ denote the components of $G - M$, where $M = \{u,v\}$. By our assumptions to date, every component in $G - M$ is a cycle component of order at least~$5$. If $k = 3$, then Dominator plays as his first vertex the vertex~$v$, thereby marking three vertices. Whatever move Staller plays as her second move, with the updated set $M$ the graph $G - M$ contains a path component of order at least~$4$, contradicting Claim~\ref{claim-3}(d). Hence, $k = 2$. Let $G_1$ be the cycle $C_{q_1}$ and let $G_2$ be the cycle $C_{q_2}$ where $5 \le q_1 \le q_2$. We note that \begin{math}
     n = 2 + q_1 + q_2,
\end{math} and so $q_1 + q_2$ is even.  Dominator now plays the vertex~$v$ on his first move, thereby marking at least two vertices. On Staller's second move, she marks at least two vertices. On Dominator's second move, he marks at least four vertices. On Staller's third move, she marks at least one vertex. Thereafter, if the game is not yet complete, then on average every played vertex marks at least two vertices. We therefore infer that \begin{math}
    \io'(G) \le 5 + \frac{1}{2}(n - 11) < \frac{1}{2}n,
\end{math}  a contradiction.~\smallqed

\medskip
By Claim~\ref{claim-9}, $k = 1$. Thus, $G - M$ is connected and is a cycle component of order~$q$ for some $q \ge 5$. Thus, $G_1 = C_q$ and $n = q+2$. Since $n$ is even, we have $q \ge 6$ and $q$ is even. Let $G_1$ be the cycle $v_1 v_2 \ldots v_q v_1$, where $vv_1$ is an edge.  Recall that by Claim~\ref{claim-2}, $\deg_G(v) \in \{2,3,4\}$, and so $v$ is adjacent to one, two or three vertices in $G_1$.

\begin{claim}
\label{claim-10}
If $q \ge 8$, then $\deg_G(v) \in \{3,4\}$. 
\end{claim}
\proof 
Suppose, to the contrary, that $q \ge 8$ and $\deg_G(v) = 2$. In this case, the vertex $v_1$ is the only vertex on the cycle $C_q$ adjacent to~$v$. Dominator can now finish the game in at most $\io(C_q)$ moves by playing the vertex~$v_1$ on his first move. Thus, $\io'(G) \le 1 + \io(C_q)$. By Theorem~\ref{thm:cycle}, if $q \bmod 5 = 0$, then \begin{math}
    \io(C_q) = \frac{2}{5}q, 
\end{math}  and so \begin{math}
    \io'(G) \le 1 + \frac{2}{5}q < \frac{1}{2}(q+2) = \frac{1}{2}n,
\end{math}  a contradiction. By Theorem~\ref{thm:cycle}, if $q \bmod 5 \ne 0$ and $q > 6$, then \begin{math}
    \io(C_q) = 2 \lceil \frac{q}{5} \rceil - 1 \le \frac{2}{5}(q+4) - 1 = \frac{1}{5}(2q+3), 
\end{math} and so \begin{math}
    \io'(G) \le 1 + \frac{1}{5}(2q+3) < \frac{1}{2}(q+2) = \frac{1}{2}n,
\end{math}  a contradiction.~\smallqed

\begin{claim}
\label{claim-11}
If $q \ge 8$, then $\deg_G(v) = 4$.
\end{claim}
\proof
Suppose, to the contrary, that $q \ge 8$ and $\deg_G(v) = 3$. In this case, $v$ is adjacent to two vertices that belong to the cycle $C_q$. Recall that $vv_1$ is an edge. Let $v_i$ be the second neighbor of $v$ that belongs to the cycle $C_q$. By symmetry, we may assume that \begin{math}
    i \in \{2,3,\ldots,\frac{q}{2},\frac{q}{2} + 1\}.
\end{math}  Dominator strategy is to play all his moves on the cycle $C_q$, starting by playing the vertex~$v_{q-2}$ as his first move. On his first move, Dominator marks three vertices, namely $v_{q-1}$, $v_{q-2}$, and $v_{q-3}$. 

If Staller does not play the vertex~$v$ on any of her remaining moves, then  Dominator can finish the game in at most $\io(C_q)$ moves, and so \begin{math}
     \io'(G) \le 1 + \io(C_q).
\end{math} In this case, as shown in the proof of Claim~\ref{claim-10}, we infer that $\io'(G) < \frac{1}{2}n$, a contradiction. Hence, Staller must play the vertex~$v$ in one of her moves.

Suppose that Staller plays the vertex~$v$ on her second move. If $q \in \{8,10\}$, then $n \in \{10,12\}$. In this case, Dominator finishes the game by playing the vertex~$v_4$ on his second move, and so \begin{math}
    \io'(G) \le 4 < \frac{1}{2}n,
\end{math}  a contradiction. Hence, $q \ge 12$. In this case, we note that Staller's second move marks at least three vertices, namely $v_1$, $v_i$ and $v_q$. After Staller's second move, with the updated set $M$ the graph $G - M$ either contains a path component of order at least~$4$, contradicting Claim~\ref{claim-3}(d) or $G-M$ is a disjoint union of two paths $P_3$ (when $q=12$ and $v_i=v_5$). In this case \begin{math}
    \io'(G) \leq 6 < \frac{1}{2}n.
\end{math} 

Suppose that Staller does not play the vertex~$v$ on her second move, implying that she plays a vertex on the cycle $C_q$. If $q = 8$, then $n = 10$, and Dominator can finish the game by playing one more move noting that $\io'(C_8) = 3$, and so in this case $\io'(G) \le 4 < \frac{1}{2}n$, a contradiction. If $q = 10$, then $n = 12$, and Dominator can finish the game by guaranteeing that at most two additional moves are played noting that $\io'(C_{10}) = 4$, and so in this case \begin{math}
    \io'(G) \le 5 < \frac{1}{2}n,
\end{math}  a contradiction. Hence, $q \ge 12$. Whatever move Staller played on the cycle as her second move, with the updated set $M$ the graph $G - M$ contains a path component of order at least~$4$, contradicting Claim~\ref{claim-3}(d).~\smallqed

\begin{claim}
\label{claim-12}
$q = 6$.
\end{claim}
\proof Suppose that $q \ge 8$. By Claim~\ref{claim-11}, we have $\deg_G(v) = 4$. If $q = 8$, then $n = 10$, and Dominator can finish the game in at most four moves, and so \begin{math}
    \io'(G) \le 4 < \frac{1}{2}n,
\end{math}  a contradiction. If $q = 10$, then $n = 12$, and Dominator can finish the game in at most five moves, and so \begin{math}
     \io'(G) \le 5 < \frac{1}{2}n,
\end{math} a contradiction. If $q = 12$, then $n = 14$ and Dominator can finish the game in at most six moves by playing, for example, the vertex~$v_1$ as his first move, and so $\io'(G) \le 6 < \frac{1}{2}n$, a contradiction. Hence, $q \ge 14$. Dominator strategy is to play all his moves on the cycle $C_q$, starting by playing the vertex~$v_1$ as his first move.

If Staller does not play the vertex~$v$ on any of her remaining moves, then  Dominator can finish the game in at most $\io(C_q)$ moves, and so $\io'(G) \le 1 + \io(C_q)$. In this case, as shown in the proof of Claim~\ref{claim-10}, we infer that $\io'(G) < \frac{1}{2}n$, a contradiction. Hence, Staller must play the vertex~$v$ in one of her moves. Dominator can guarantee that the game finishes in at most~$2 + \io'(C_{q})$ moves. Thus, \begin{math}
    \io'(G) \leq 2 + \io'(C_{q}). 
\end{math}  Recall that $q \ge 14$ and $q$ is even.

If $q = 14$, then by Theorem~\ref{thm:cycle}, we have \begin{math}
    \io(C_q) = 2 \lceil \frac{q}{5} \rceil - 1 = \frac{2}{5}(q+1) - 1 = \frac{1}{5}(2q-4), 
\end{math} and so \begin{math}
    \io'(G) \le 2 + \frac{1}{5}(2q-4) < \frac{1}{2}(q+2) = \frac{1}{2}n,
\end{math}  a contradiction. Moreover if $q = 16$, then Dominator strategy of playing all his moves on the cycle $C_q$, starting by playing the vertex~$v_1$ as his first move, allows him to complete the game in at most eight moves, and so \begin{math}
    \io'(G) \le 8 < \frac{1}{2}n,
\end{math}  a contradiction. Hence, $q \ge 18$.

By Theorem~\ref{thm:cycle}, if $q \bmod 5 = 0$ (and $q > 10$ is even), then $\io(C_q) = \frac{2}{5}q$, and so \begin{math}
     \io'(G) \le 2 + \frac{2}{5}q < \frac{1}{2}(q+2) = \frac{1}{2}n,
\end{math} a contradiction. 

By Theorem~\ref{thm:cycle}, if $q \bmod 5 \ne 0$ and $q \ge 18$, then \begin{math}
    \io(C_q) = 2 \lceil \frac{q}{5} \rceil - 1 \le \frac{2}{5}(q+4) - 1 = \frac{1}{5}(2q+3),
\end{math}  and so \begin{math}
    \io'(G) \le 2 + \frac{1}{5}(2q+3) < \frac{1}{2}(q+2) = \frac{1}{2}n,
\end{math}  a contradiction. Hence the case $q \ge 8$ cannot occur, implying that $q = 6$.~\smallqed

\medskip
We now return to the proof of Theorem~\ref{thm:1/2bound-Sgame} one final time. By Claim~\ref{claim-12}, we have $q = 6$ and $n = 8$. Noting that $\io'(G) = \frac{1}{2}n = 4$, we readily infer that \begin{math}
    G \in \{F_7,F_8,F_9,F_{10},F_{11}\}.
\end{math}  This completes the proof of Theorem~\ref{thm:1/2bound-Sgame}.
     \end{proof}

\medskip
\begin{remark}
\label{rem:connected-not-needed}    
The inequalities of Theorems~\ref{thm:1/2bound} and~\ref{thm:1/2bound-Sgame} hold true for arbitrary graphs, that is, they need not to be connected. This can be checked by observing that in the proofs of the inequalities, connectivity is not an issue. \end{remark}
 
\section{The isolation game played in trees}
\label{sec:tree-bound}

As a consequence of Theorem~\ref{thm:1/2bound}, we have the following upper bound on the game isolation number of a tree.

\begin{corollary}
\label{cor:tree-bound-1}
If $T$ is a tree of order~$n$, then $\io(T) \le \frac{1}{2}n$, with equality if and only if $T = K_2$.  
\end{corollary}

In this section, we improve the $\frac{1}{2}n$-upper bound in Corollary~\ref{cor:tree-bound-1} to a $\frac{5}{11}n$-upper bound when $T$ is a tree of order~$n \ge 3$. For this purpose, we use a proof technique from~\cite{GoHe-25} used to establish upper bounds on the isolation number of a graph. The technique uses  an \emph{isolation residual graph} defined as follows. 
 	
Given a graph $G$ and a subset $S \subseteq V(G)$, the \emph{isolation residual graph} denoted by $\residual$ is obtained as follows. Its vertex set is $V(G)$. Further:
\begin{enumerate}
\item
A vertex in $G-N[S]$ with a neighbor in $G-N[S]$  is colored \emph{white}.  
\item
A vertex in $N[S]$ with a white neighbor is colored \emph{blue}.
\item
All other vertices are colored \emph{red}.
\item
The edge set of $\residual$ consists of all edges of $G$ that are incident with at least one white vertex.
\end{enumerate}

As remarked in~\cite{GoHe-25}, the white vertices still need to be isolated, the blue vertices might be useful to do this isolating, and the red vertices are no longer relevant. We denote the set of red vertices by~$R$, the set of blue vertices by $B$, and the set of white vertices by $W$. We observe that these sets form a (weak) partition of $V(G)$. We shall need the following elementary properties of a residual graph. 

\begin{observation}
\label{ob-residual}
The following holds in the isolation residual graph $\residual$.
\begin{enumerate}
\item
Every vertex of $S$ is colored red.
\item
Every red vertex is isolated, but no white nor blue vertex is isolated.
\item
The degree of a white vertex in $\residual$ equals its degree in $G$.
\item
The set $S$ is an isolating set in $G$ if and only if there is no white vertex.
\end{enumerate}
\end{observation}

We are now in a position to prove the following result on the game isolation number of a tree.  

\begin{theorem}
\label{thm:tree-bound-2}
If $T$ is a tree of order~$n \ge 3$, then \begin{math}
    \io(T) \le \frac{5}{11}n.
\end{math} 
\end{theorem}
\begin{proof}
For a set $S \subseteq V(T)$, consider the isolation residual graph $\resT$ and the sets of white, blue, and red vertices in $\resT$. We define the weight function 
\begin{displaymath}
    \w(\resT)= 5 |W| +3|B|. 
\end{displaymath} 
Equivalently, we associate every white, blue, and red vertex with a weight of $5$, $3$, and $0$, respectively, and then consider the sum of the weights in $\residual$.

When the game starts, $S=\emptyset$ and $\w(T_\emptyset^\iota)=5n$, while $\w(\resT)=0$ holds  when the game is over. Further, if $v$ is a playable vertex in the isolation residual graph $\resT$, then \begin{math}
    \w(\resT) \ge \w(T_{S \cup \{v\}}^\iota)
\end{math} and we use the notation \begin{math}
    \xi(v) = \w(\resT)- \w(T_{S \cup \{v\}}). 
\end{math} We will show that Dominator can ensure that $\w(\resT)$ decreases by at least $22$ in every two consecutive moves of Dominator and Staller. To prove this, we consider three stages in the D-game. In the graph $\resT[W]$, induced by the white vertices of $\resT$, a \emph{big component} is a component of order at least~$4$. If a component is not big, it is called \emph{small}. If $T$ is a tree, a small component of $\resT[W]$ is either a $P_2$- or a $P_3$-subgraph.  

\paragraph{Stage~$1$.} A move of Dominator belongs to Stage~$1$ if $\resT[W]$ contains a big component before this move. A move of Staller belongs to Stage~$1$, if so does Dominator's previous move.

We first state that every legal move decreases $\w(\resT)$ by at least $5$. Indeed, if a white vertex $v$ is played, it becomes red and $\xi(v) \ge 5$. If a blue vertex $u$ is played that has a white neighbor $u'$, then $u$ becomes red and $u'$ becomes blue or red, which gives \begin{math}
    \xi(u) \ge 3 + (5-3)=5.
\end{math}  In particular, each move of Staller decreases the weight of $\resT$ by at least~$5$.

Suppose that Dominator plays in $\resT$ in Stage~$1$. The induced subgraph $\resT[W]$ then contains a component $T_i$ on at least four vertices. If $T_i$ is a star $K_{1,k}$ for a $k \ge 3$, then Dominator chooses the center $v$. Since all vertices in the star become red, \begin{math}
    \xi(v) \ge 5(k+1) \ge 20.
\end{math}  If $T_i$ is not a star, it contains a diametrical path $v_1\dots v_d$ with $d \ge 4$. Then $P\colon v_1v_2v_3v_4$ is a $4$-path such that $v_1$ is a leaf in $T_i$ and if $v_2$ has a neighbor $v' \notin V(P)$, then $v'$ is a leaf in $T_i$.  
Dominator now plays $v_3$. With this move $v_1$, $v_2$, $v_3$ (more generally, $N_{T_i}[v_2]$) become red and $v_4$ is colored blue or red. Consequently, \begin{math}
    \xi(v_3) \ge 3\times 5+ (5-3)= 17.
\end{math}  This shows that Dominator can always decrease $\w(\resT)$ by at least~$17$, when the move belongs to Stage~$1$. We may conclude that Dominator has a strategy to ensure that the weight decreases by at least $17+5=22$ in every two consecutive moves in Stage~$1$.
\medskip

Before discussing Stage~$2$, we introduce a new vertex type in the isolation residual graph $\resT$. We say that a vertex $u$ is \emph{light blue}, and belongs to $B^\ell$, if it became blue during Stage~$2$. Light blue vertices are associated with a weight of $1$; that is, the new weight function we use in Stages~$2$ and $3$ is
\begin{displaymath}
    \w(\resT)= 5 |W| +3|B\setminus B^\ell| + |B^\ell|. 
\end{displaymath} 
Note that all light blue vertices are of degree~$1$. Further, we say that a component $T_i$ of $\resT[W]$ is \emph{special}, if it is a $P_2$-component and each adjacent blue (or light blue) vertex is of degree~$1$ in $\resT$. 

\paragraph{Stage 2.} A move of Dominator belongs to Stage~$2$, if the subgraph $\resT[W]$ in $\resT$ before this move consists only of small components and not all of them are special. A move of Staller belongs to Stage~$2$, if so does Dominator's previous move.

We first prove that every move $v$ of Staller decreases $\w(\resT)$ by at least~$7$ in Stage~$2$. If $v$ is a white vertex, then an entire $P_2$- or $P_3$-component of $G[W]$ becomes red and \begin{math}
    \xi(v) \ge 2\times 5=10.
\end{math}  The same is true if $v \in B$ and a white neighbor $v'$ belongs to a $P_2$-component. Assume $v \in B$ and that a neighbor $v'$ of $v$ is from a $P_3$-component of $\resT[W]$. We remark that in this case $v \in B\setminus B^\ell$ because there were no components in $\resT[W]$ on more than three vertices at the beginning of Stage~$2$, and so $v$ was not white at that time. 
If $v'$ is the central vertex of the $P_3$-component, then $v'$ and the entire $P_3$-component becomes red with this move and we have \begin{math}
    \xi(v) \ge 3+ 3 \times 5=18.
\end{math}  If $v'$ is a leaf in the $P_3$-component, then $v$ becomes red and $v'$ becomes light blue with this move. Therefore, \begin{math}
    \xi(v) \ge 3+(5-1)=7,
\end{math}  completing the proof of the statement.

Consider now Dominator's move. If a $P_3$-component is present in $\resT[W]$, he can choose a vertex $v$ from such a component $T_i$. All three white vertices of $T_i$ become red, which results in \begin{math}
    \xi(v) \ge 3 \times 5=15.
\end{math}  If there is no $P_3$-component, then a non-special $P_2$-component exists. Dominator may choose a blue vertex $u$ that has neighbors $u_1$ and $u_2$ from two different $P_2$-components. Playing $u$, both components and $u$ itself become red and \begin{math}
    \xi(u) \ge 4\times 5+3=23.
\end{math}  We may conclude that Dominator is always able to decrease $\w(\resT)$ by at least~$15$ with a move in Stage~$2$ and therefore, every two consecutive moves together decrease the weight of the isolation residual graph by at least $7+15=22$.

\paragraph{Stage~$3$.} Stage~$3$ starts when all components of $\resT[W]$ are special and ends when the game is over.

Each move of a player removes a $P_2$-component from $\resT$ by recoloring both vertices with red. Further, by Observation~\ref{ob-residual}(c) and  the condition $n \ge 3$ in the theorem, each special component is adjacent to a blue vertex of degree~$1$, which also becomes red together with the component. We infer that $\w(\resT)$ decreases by at least $2\times 5+1=11$ with each move played in Stage~$3$, no matter whether the move is played by Dominator or by Staller.
\medskip

By the analysis of the stages of the game when Dominator follows the described strategy, we obtained that, on average, each move decreases $\w(\resT)$ by at least $11$. Let $t$ denote the number of vertices played in the game. As the game started with $\w(T_\emptyset^\iota)=5n$ and finished with $\w(\resT)=0$, we conclude $11t \leq 5n$ and then, \begin{math}
    \io(T) \leq t \leq \frac{5}{11} n
\end{math}  follows.     
\end{proof}

\section{On the 3/7-bound}
\label{sec:3/7}

\cite{bdj-2024} posed the following conjecture.

\begin{conjecture}\label{con:3/7}
    For any graph $G$ with no $K_2$-components, \begin{displaymath}
        \io(G) \leq \left\lceil\frac{3}{7}n(G)\right\rceil \text{ and } \io'(G) \leq \left\lceil\frac{3}{7}n(G)\right\rceil.
    \end{displaymath}
\end{conjecture}

If the conjecture holds true, then it is sharp as the following example from~\cite{bdj-2024} demonstrates. Start with an arbitrary graph $G$ and attach to each vertex of it two disjoint paths of length three. Besides this family of graphs and graphs similar to them (say, obtained from such a graph by removing a pendant vertex), the authors could not find any other infinite family of graphs whose game isolation number is close to $3n(G)/7$. Here, we present another such family. 

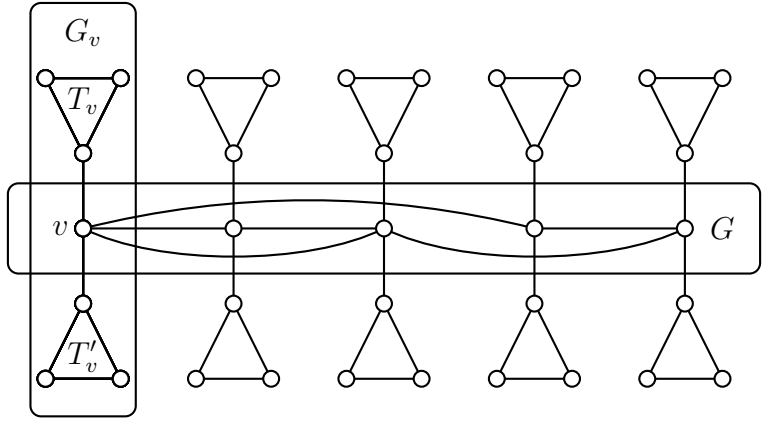
\begin{figure}[ht!]
\begin{center}
\begin{tikzpicture}[scale=1.0,style=thick,x=1cm,y=1cm]
\def\vr{3pt}
\begin{scope}[xshift=-0cm, yshift=0cm] 
\coordinate(x1) at (0.0,0.0);
\coordinate(y1) at (0.0,1.0);
\coordinate(y2) at (-0.5,2.0);
\coordinate(y3) at (0.5,2.0);
\coordinate(z1) at (0.0,-1.0);
\coordinate(z2) at (-0.5,-2.0);
\coordinate(z3) at (0.5,-2.0);
\draw (x1) -- (y1) -- (y2) -- (y3) -- (y1);  
\draw (x1) -- (z1) -- (z2) -- (z3) -- (z1);  
\draw (0,0) -- (2,0);
\draw (2,0) -- (4,0);
\draw (6,0) -- (8,0);
\draw (0,0) .. controls (2,0.5) and (4,0.5) .. (6,0);
\draw (4,0) .. controls (5,-0.5) and (7,-0.5) .. (8,0);
\draw (0,0) .. controls (1,-0.5) and (3,-0.5) .. (4,0);
\draw(x1)[fill=white] circle(\vr);
\foreach \i in {1,2,3} 
{
\draw(y\i)[fill=white] circle(\vr);
\draw(z\i)[fill=white] circle(\vr);
}
\node at (8.5,0.0) {$G$};
\node at (-0.3,0.0) {$v$};
\node at (0.0,1.7) {$T_v$};
\node at (0.0,-1.7) {$T_v'$};
\node at (0.0,2.6) {$G_v$};
\draw[rounded corners] (-1.0, -0.6) rectangle (9.0, 0.6);
\draw[rounded corners] (-0.7, -2.5) rectangle (0.7, 3.0);

\end{scope}
\begin{scope}[xshift=2cm, yshift=0cm] 
\coordinate(x1) at (0.0,0.0);
\coordinate(y1) at (0.0,1.0);
\coordinate(y2) at (-0.5,2.0);
\coordinate(y3) at (0.5,2.0);
\coordinate(z1) at (0.0,-1.0);
\coordinate(z2) at (-0.5,-2.0);
\coordinate(z3) at (0.5,-2.0);
\draw (x1) -- (y1) -- (y2) -- (y3) -- (y1);  
\draw (x1) -- (z1) -- (z2) -- (z3) -- (z1);  
\draw(x1)[fill=white] circle(\vr);
\foreach \i in {1,2,3} 
{
\draw(y\i)[fill=white] circle(\vr);
\draw(z\i)[fill=white] circle(\vr);
}
\end{scope}
\begin{scope}[xshift=4cm, yshift=0cm] 
\coordinate(x1) at (0.0,0.0);
\coordinate(y1) at (0.0,1.0);
\coordinate(y2) at (-0.5,2.0);
\coordinate(y3) at (0.5,2.0);
\coordinate(z1) at (0.0,-1.0);
\coordinate(z2) at (-0.5,-2.0);
\coordinate(z3) at (0.5,-2.0);
\draw (x1) -- (y1) -- (y2) -- (y3) -- (y1);  
\draw (x1) -- (z1) -- (z2) -- (z3) -- (z1);  
\draw(x1)[fill=white] circle(\vr);
\foreach \i in {1,2,3} 
{
\draw(y\i)[fill=white] circle(\vr);
\draw(z\i)[fill=white] circle(\vr);
}
\end{scope}
\begin{scope}[xshift=6cm, yshift=0cm] 
\coordinate(x1) at (0.0,0.0);
\coordinate(y1) at (0.0,1.0);
\coordinate(y2) at (-0.5,2.0);
\coordinate(y3) at (0.5,2.0);
\coordinate(z1) at (0.0,-1.0);
\coordinate(z2) at (-0.5,-2.0);
\coordinate(z3) at (0.5,-2.0);
\draw (x1) -- (y1) -- (y2) -- (y3) -- (y1);  
\draw (x1) -- (z1) -- (z2) -- (z3) -- (z1);  
\draw(x1)[fill=white] circle(\vr);
\foreach \i in {1,2,3} 
{
\draw(y\i)[fill=white] circle(\vr);
\draw(z\i)[fill=white] circle(\vr);
}
\end{scope}
\begin{scope}[xshift=8cm, yshift=0cm] 
\coordinate(x1) at (0.0,0.0);
\coordinate(y1) at (0.0,1.0);
\coordinate(y2) at (-0.5,2.0);
\coordinate(y3) at (0.5,2.0);
\coordinate(z1) at (0.0,-1.0);
\coordinate(z2) at (-0.5,-2.0);
\coordinate(z3) at (0.5,-2.0);
\draw (x1) -- (y1) -- (y2) -- (y3) -- (y1);  
\draw (x1) -- (z1) -- (z2) -- (z3) -- (z1);  
\draw(x1)[fill=white] circle(\vr);
\foreach \i in {1,2,3} 
{
\draw(y\i)[fill=white] circle(\vr);
\draw(z\i)[fill=white] circle(\vr);
}
\end{scope}
\begin{scope}[xshift=-0cm, yshift=0cm] 
\coordinate(x1) at (0.0,0.0);
\coordinate(y1) at (0.0,1.0);
\coordinate(y2) at (-0.5,2.0);
\coordinate(y3) at (0.5,2.0);
\coordinate(z1) at (0.0,-1.0);
\coordinate(z2) at (-0.5,-2.0);
\coordinate(z3) at (0.5,-2.0);
\draw (x1) -- (y1) -- (y2) -- (y3) -- (y1);  
\draw (x1) -- (z1) -- (z2) -- (z3) -- (z1);  
\draw(x1)[fill=white] circle(\vr);
\foreach \i in {1,2,3} 
{
\draw(y\i)[fill=white] circle(\vr);
\draw(z\i)[fill=white] circle(\vr);
}
\end{scope}
\begin{scope}[xshift=-0cm, yshift=0cm] 
\coordinate(x1) at (0.0,0.0);
\coordinate(y1) at (0.0,1.0);
\coordinate(y2) at (-0.5,2.0);
\coordinate(y3) at (0.5,2.0);
\coordinate(z1) at (0.0,-1.0);
\coordinate(z2) at (-0.5,-2.0);
\coordinate(z3) at (0.5,-2.0);
\draw (x1) -- (y1) -- (y2) -- (y3) -- (y1);  
\draw (x1) -- (z1) -- (z2) -- (z3) -- (z1);  
\draw(x1)[fill=white] circle(\vr);
\foreach \i in {1,2,3} 
{
\draw(y\i)[fill=white] circle(\vr);
\draw(z\i)[fill=white] circle(\vr);
}
\end{scope}
\begin{scope}[xshift=-0cm, yshift=0cm] 
\coordinate(x1) at (0.0,0.0);
\coordinate(y1) at (0.0,1.0);
\coordinate(y2) at (-0.5,2.0);
\coordinate(y3) at (0.5,2.0);
\coordinate(z1) at (0.0,-1.0);
\coordinate(z2) at (-0.5,-2.0);
\coordinate(z3) at (0.5,-2.0);
\draw (x1) -- (y1) -- (y2) -- (y3) -- (y1);  
\draw (x1) -- (z1) -- (z2) -- (z3) -- (z1);  
\draw(x1)[fill=white] circle(\vr);
\foreach \i in {1,2,3} 
{
\draw(y\i)[fill=white] circle(\vr);
\draw(z\i)[fill=white] circle(\vr);
}
\end{scope}
\begin{scope}[xshift=-0cm, yshift=0cm] 
\coordinate(x1) at (0.0,0.0);
\coordinate(y1) at (0.0,1.0);
\coordinate(y2) at (-0.5,2.0);
\coordinate(y3) at (0.5,2.0);
\coordinate(z1) at (0.0,-1.0);
\coordinate(z2) at (-0.5,-2.0);
\coordinate(z3) at (0.5,-2.0);
\draw (x1) -- (y1) -- (y2) -- (y3) -- (y1);  
\draw (x1) -- (z1) -- (z2) -- (z3) -- (z1);  
\draw(x1)[fill=white] circle(\vr);
\foreach \i in {1,2,3} 
{
\draw(y\i)[fill=white] circle(\vr);
\draw(z\i)[fill=white] circle(\vr);
}
\end{scope}
\begin{scope}[xshift=-0cm, yshift=0cm] 
\coordinate(x1) at (0.0,0.0);
\coordinate(y1) at (0.0,1.0);
\coordinate(y2) at (-0.5,2.0);
\coordinate(y3) at (0.5,2.0);
\coordinate(z1) at (0.0,-1.0);
\coordinate(z2) at (-0.5,-2.0);
\coordinate(z3) at (0.5,-2.0);
\draw (x1) -- (y1) -- (y2) -- (y3) -- (y1);  
\draw (x1) -- (z1) -- (z2) -- (z3) -- (z1);  
\draw(x1)[fill=white] circle(\vr);
\foreach \i in {1,2,3} 
{
\draw(y\i)[fill=white] circle(\vr);
\draw(z\i)[fill=white] circle(\vr);
}
\end{scope}

\end{tikzpicture}
\caption{Graph $\widehat{G}$}
\label{fig:3/7-family}
\end{center}
\end{figure}

If $G$ is a graph, then let $\widehat{G}$ be a graph defined as follows.  For any vertex $v\in V(G)$, add two disjoint triangles $T_v$ and $T_v'$, and add an edge between $v$ and a vertex of $T_v$ and an edge between $v$ and a vertex of $T_v'$, see Fig.~\ref{fig:3/7-family}.

\begin{proposition}
\label{prop:G-hat}
If $G$ is a graph, then 
\begin{displaymath}
    \io(\widehat{G}) = \io'(\widehat{G}) = \frac{3}{7} n(\widehat{G})\,.
\end{displaymath} 
\end{proposition}

\begin{proof}
 For a vertex $v\in V(G)$, let $G_v$ be the subgraph of $\widehat{G}$ induced by \begin{math}
     \{v\}\cup V(T_v)\cup V(T_v'),
 \end{math}  see Fig.~\ref{fig:3/7-family} again. 

Consider the following strategy of Dominator that guarantees that not more than three vertices are played in each $G_v$. Assume that the last move of Staller is a vertex from $G_v$, for some $v\in V(G)$. If there are vertices of $G_v$ which remain playable, then Dominator plays a vertex of $G_v$ such that no more than three vertices from $G_v$ will be played when the game finishes.  He can achieve this goal since if Staller played $v$, then Dominator selects a vertex from $T_v$ or a vertex from $T_v'$, and if Staller selects some vertex from $V(T_v)\cup V(T_v')$, then Dominator replies by playing $v$ if possible, otherwise a vertex from the triangle where no vertex has been played yet. Note that Dominator can also follow this strategy for his first move in the D-game, that is, he selects some vertex $v\in V(G)\subseteq V(\widehat{G})$. On the other hand, if no vertex of $G_v$ remains playable after the last move of Staller, Dominator applies analogous strategy in some $G_{u}$, $u\ne v$, where there still exist playable vertices. Note that this strategy of Dominator can be applied no matter whether the D-game or the S-game is played. We can conclude that \begin{math}
    \io(\widehat{G})\le \frac{3}{7} n(\widehat{G})
\end{math}  as well as \begin{math}
  \io'(\widehat{G}) \le \frac{3}{7} n(\widehat{G}).   
\end{math} 

Consider now the following strategy of Staller. Her goal is that in each $G_v$, the vertex $v$ is selected during the game, either by Dominator or by Staller. So assume that the last move of Dominator is a vertex from $G_v$, for some $v\in V(G)$, and that no other vertex of $G_v$ has yet been played. If  Dominator selected $v$, then Staller replies by selecting a neighbor of $v$ in $T_v$. On the other hand, if Dominator selected a vertex different from $v$, then Staller selects $v$. In either case, this strategy of Staller guarantees that by the end of the game, three vertices from $G_v$ are selected. In addition, if Staller is forced to play first on some $G_v$, then she selects $v$, thus forcing again that three vertices from $G_v$ will be played. We can conclude that \begin{math}
    \io(\widehat{G})\ge \frac{3}{7} n(\widehat{G})
\end{math}  as well as \begin{math}
    \io'(\widehat{G}) \ge \frac{3}{7} n(\widehat{G}),
\end{math}  and we are done. 
  \end{proof}

Both the previously known extremal family (based on attaching two paths of length three) and our construction (based on attaching two triangles) share a common structural feature. Each vertex of the original graph is extended into a local configuration on seven vertices that requires at least three vertices in any optimal isolating set. Consequently, these constructions can be viewed as decomposing the graph into disjoint 7-vertex gadgets contributing exactly 3 to the game isolation number, which explains the emergence of the $\frac{3}{7}$ ratio. This suggests that the conjectured upper bound is driven by such local configurations, rather than global properties of the graph, and that extremal graphs may be characterized by the presence of these forcing gadgets.

\medskip
We conclude the paper with the following slight refinements of Conjecture~\ref{con:3/7}.

\begin{conjecture}
\label{con:first}
If $G$ is a graph of order $n$ with $\delta(G) \ge 2$, then $\io(G) \le \frac{3}{7}n$ and $\io'(G) \le \frac{3}{7}n$.  
\end{conjecture}

If Conjecture~\ref{con:first} is correct, then we have seen above that it is tight and achievable for connected graphs $G$ of arbitrarily large order. 

\begin{conjecture}
\label{con:second}
If $T$ is a tree of order $n \ge 3$, then $\io(T) \le \frac{3}{7}n$.  
\end{conjecture}

\acknowledgements
\label{sec:ack}
A lot of work on this paper has been done during the Workshop on Games on Graphs III, June 2025, Rogla, Slovenia, the authors thank the Institute of Mathematics, Physics and Mechanics, Ljubljana, Slovenia for supporting the workshop. 

\nocite{*}
\bibliographystyle{abbrvnat}
\bibliography{bujtas-dmtcs}
\label{sec:biblio}

\end{document}